\input amstex
\input amssym
\documentstyle{amsppt}
\def\e{\epsilon}
\def\a{\alpha}
\def\b{\beta}
\def\th{\theta}
\def\d{\delta}
\def\G{\Gamma}
\def\g{\gamma}

\def\D{\Delta}

\def\s{\sigma}

\def\x{\times}

\def\f{\flushpar}

\def\om{\omega}
\def\Om{\Omega}
\def\B{\Cal B}

\def\mathcal{\Cal}
\def\({\biggl(}
\def\){\biggr)}

\def\<{\bold\langle}
\def\>{\bold\rangle}

\def \r{\Cal R}\def\bdy{\partial}

\def\bul{\smallskip\f$\bullet\ \ \ $}
\document\topmatter
\title Exchangeable, Gibbs and equilibrium measures for
 Markov subshifts
  \endtitle
\author J. Aaronson,\   H. Nakada \endauthor
 \address[Jon. Aaronson]{\ \ School of Math. Sciences, Tel Aviv University,
69978 Tel Aviv, Israel.}
\endaddress
\email{aaro\@tau.ac.il}\endemail\address[Hitoshi Nakada]{\ \ Dept.
of Math., Keio University,Hiyoshi 3-14-1 Kohoku,
 Yokohama 223, Japan}\endaddress
\email{nakada\@math.keio.ac.jp}\endemail\abstract We study a class
of  strongly irreducible,
 multidimensional, topological Markov shifts, comparing two notions of ``symmetric measure":
 exchangeability and the Gibbs (or conformal) property. We show
   that equilibrium measures for such shifts (unique and weak Bernoulli in the one dimensional case)
    exhibit a variety of   spectral properties.
\endabstract\thanks\copyright 2005. 
\endthanks\subjclass 60G09,\ 37D35 (37A15,\ 37A20,\ 37A40,\ 60G60,\ 60J10)\endsubjclass
 \endtopmatter
\rightheadtext{Exchangeability}\heading { Introduction}\endheading

  Let $S$ be a finite set of {\it spins} and let $\Gamma$ be a countable set of {\it sites}.
The  {\it tail} (or {\it Gibbs} or {\it homoclinic}) relation on a
{\it configuration set} (or {\it lattice system}) $X\in\B(S^\G)$ (the Borel subsets of $X$)
is defined by
$$\goth T(X):=\{(x,y)\in X\x X:\ \exists\ F\subset\G,\ \# F<\infty,\ x_{F^c}=y_{F^c}\}.$$
Here, for $x\in S^\G,\ \Lambda\subset\G$, $x_\Lambda\in S^\Lambda$
is the {\it $\Lambda$-restriction} of $x$ to $\Lambda$, defined by
$(x_\Lambda)_j=x_j\ \ \ (j\in\Lambda)$ (and the collection of
$\Lambda$-restrictions is $X_\Lambda:=\{x_\Lambda:\ x\in
X\}\subset S^\Lambda$).

  The {\it exchangeable relation} on $X$ is
  $$\align \Cal E(X)& :=\{(x,y)\in X\x X:\ \exists\ F\subset\G,\ \# F<\infty,\
  x_{F^c}=y_{F^c},\\ & \ \ \ \ \ \ \ \ \ \ \ \ \exists\ \text{ a permutation }\s:F\to F,\ y_i=x_{\s(i)}\ \forall\ i\in F\}\\ &
   \subseteq\goth T(X).\endalign$$
Both the exchangeable and tail relations are  countable
equivalence relations in the sense of \cite{FM}  (Borel measurable
equivalence relations with countable equivalence classes).
\par In this paper, we consider $\mathcal
E(X)$- and $\goth T(X)$-invariant measures when $\G=\Bbb Z^d$ and
$X$ is a {\it $\Bbb Z^d$-topological Markov  shift}.

The {\it shift action}  $T$ on $S^{\Bbb Z^d}$  is the $\Bbb Z^d$
 action $T_k:S^{\Bbb Z^d}\to S^{\Bbb Z^d}\ \ (k\in\Bbb Z^d)$
 defined by $(T_kx)_n:=x_{n+k}\ \ \ (k,n\in\Bbb Z^d)$. A $\Bbb Z^d$-{\it
 subshift} is a closed subset  $X\subset S^{\Bbb Z^d}$ which is  $T$-invariant
 ($T_kX=X\ \forall\ k\in\Bbb Z^d$).

 In order to define a $\Bbb Z^d$-topological Markov  shift, consider
 $\Bbb Z^d$  equipped with the norm
 $\|n\|=\|(n_1,\dots,n_d)\|:=\max_{1\le k\le d}|n_k|$ and let
$B(n,r):=\{k\in\Bbb Z^d:\ \|k-n\|\le
r\}=\prod_{k=1}^d[n_k-r,n_k+r]\cap\Bbb Z^d.$ Given a set
$F\subset\Bbb Z^d$, let $F^o:=\{x\in F:\ B(x,1)\subset F\},\ \bdy
F:=F\setminus F^o$.

 A $\Bbb Z^d$-{\it topological Markov  shift} (TMS) is a subshift $X\subset S^{\Bbb Z^d}$
determined by {\it nearest neighbor constraints}: there
is a subset $A\subset S\x S^{\bdy B(0,1)}$ so
 that
$$X=\{x\in S^{\Bbb Z^d}:\ (x_k,x_{k+\bdy B(0,1)})\in A\ \forall\
k\in\Bbb Z^d\}.$$

We consider two kinds of ``naturally symmetric measure" on $X$: an
\bul {\it exchangeable measure } being $\Cal E(X)$-invariant, and
a \bul {\it local Gibbs measure} (or {\it state}): a conformal
measure with locally determined potential (see \cite{R} and below).

Site-Gibbs measures (where the potential  is determined by the
site) are automatically exchangeable.

\par By the   De Finetti-Hewitt-Savage theorem (\cite{He-Sa}),  the
$\Cal E(S^{\Bbb Z^d })$-invariant, ergodic probabilities on  $S^{
\Bbb Z^d }$ are precisely the  stationary product measures (which
are the extremal, site-Gibbs measures on  $S^{ \Bbb Z^d }$).

 Theorems 1 and  2 show that under certain conditions,    an  ergodic
exchangeable measure on a $\Bbb Z^d$-TMS which is {\it global}
 (i.e. globally supported in the sense that every open set has positive measure)
 is a site-Gibbs measure.
\par The notion of "restricted exchangeability" appears in
\cite{PS} where the exchangeable, shift invariant, ergodic
probabilities for $\Bbb Z$-TMS's were identified and an extensive
bibliography on the subject is provided. Exchangeable measures on
one sided TMS's were considered in \cite{ANSS1}. For previous
results concerning the multidimensional subshift cases see
\cite{S2}.

The simplest kind of site-Gibbs measure on $X$ is a $\goth
T(X)$-invariant measure. These exist and are Markov with uniform
specifications. They are unique when $d=1$, but not when $d\ge 2$
(see \cite{BS1}, \cite{BS2}, \cite{BS3}) .\par A $T$-invariant,
$T$-ergodic, $\goth T(X)$-invariant  measure on a strongly
irreducible $\Bbb Z^d$-TMS is an {\it equilibrium measure} having
maximal entropy. Theorem 3 (being a reconsideration of a
Burton-Steif construction) shows that- \bul the equilibrium
measures of strongly irreducible $\Bbb Z^3$-TMS's  \bul the Gibbs
measures on $\{-1,1\}^{\Bbb Z^3}$ with a nearest neighbour
potential

\smallskip
exhibit a variety of spectral properties.

\smallskip\par The main methods of this paper are the
 theories of cocycles and equivalence relations as
introduced in \cite{S1} and  \cite{FM} (respectively). Definitions
can be found on \S0 (after this introduction). The theorems are
stated in \S1 and the rest of the paper is devoted to proofs and
examples. \heading \S0 Definitions\endheading \subheading{0.1
Countable equivalence relations}  As in \cite{FM}, a {\it
countable equivalence relation} on a Polish space $X$ is a
subset $\r\in\B(X\x X)$ which is an equivalence relation with
countable equivalence classes.

Both the exchangeable and tail relations are  countable
equivalence relations on $X\in\B(S^\G)$ (where $S$ is finite and
$\G$ is countable) .

If  $G$ is a countable group of measurable, invertible
transformations of $X$, then
$$\r_G:=\{(x,Rx):\ R\in G,\ x\in X\}$$ is a countable equivalence
relation on $X$.

By [FM], every countable, equivalence relation is of form $\r_G$
for some $G$.

In this paper, we consider various dynamical properties of
countable equivalence relations. Unless stated otherwise, by a
dynamical property of a countable group $G$ of measurable,
invertible transformations we mean the corresponding property of
$\r_G$.
\par Let $\r$ be  a countable equivalence relation on $X$. A
 {\it $\r$-holonomy} is a Borel isomorphism
$\Phi:B\to C,\ \ (B,\ C\in\mathcal B)$ with $ (x,\Phi(x))\in \r\
\forall\  x\in B$ (shorthand: $\Phi:B\overset{\r}\to\rightarrow
C$). A  $\r$-holonomy $\Phi:B\overset{\r}\to\rightarrow C$ is
called {\it topological } if $B,\ C\subset X$ are open and
$\Phi:B\to C$ is a homeomorphism.

A collection $\goth C$ of $\mathcal R$- holonomies {\it generates
$\mathcal R$} if for each $(x,y)\in\mathcal R,\ \exists\
\Phi\in\goth C,\ \Phi:B\overset{\r}\to\rightarrow C,\ x\in B,\
\Phi(x)=y.$ Any countable equivalence relation on $X$ is generated
by a countable collection of $\mathcal R$- holonomies (as shown in
\cite{FM}). The collection of $\r$-holonomies is denoted by
$[[\r]]$ and known as the {\it groupoid} of $\r$. The {\it full
group} of $\r$ is collection of globally defined $\r$-holonomies:
$$[\r]:=\{\Phi\in [[\r]]:\ \Phi:X\overset{\r}\to\rightarrow X\}.$$
It is a group under composition.

A countable equivalence relation $\r$ on a topological space $X$
is called a {\it topological equivalence relation on $X$} if it is
generated by a countable set  of topological
 $\r$-holonomies.
\par In case $X\subset S^\G$, ($\G$ countable) we consider {\it topological
cylinder holonomies}. For $F\subset\G$ finite, and $a\in X_F$, we
define the {\it $F$-cylinder}   (with configuration $a$) as
$$[a]_F:=\{x\in X:\ x_F=a\}.$$ This is a clopen subset of $X$.

We call the $F$-configurations  $a,b\in X_F$ {\it compatible} if
$\forall\ x\in [a]_F,\ \exists\ y\in [b]_F$ with $y_{F^c}=x_{F^c}$
and vice versa. \par The pair $a,b\in X_F$ is compatible iff the
map $\pi:(x_{F^c},a)\mapsto (x_{F^c},b)$ is a homeomorphism
$\pi:[a]_F\to [b]_F$, i.e. $\pi:[a]_F\overset{\goth
T(X)}\to\longrightarrow  [b]_F$ is a topological holonomy. Call
such holonomies
 {\it topological
cylinder holonomies}. If $\goth T(X)$ is generated by topological
cylinder holonomies, then it is  a topological equivalence
relation.

Suppose that $X\subset
 S^{\Bbb Z^d}$ ($d\ge 1$) is a $\Bbb Z^d$-TMS. Let $F\subset\Bbb Z^d$ be finite and suppose that
 $a,b\in X_F$. If $a_{\bdy F}=b_{\bdy F}$, then $a, b$ are compatible and are connected by a topological cylinder holonomy.
  $\goth T(X)$ is generated by such holonomies and is thus  a topological equivalence relation on
 $X$.

\subheading{0.2 Measures}
\par We denote by $\Cal P(X)$ the collection of  probability measures on $X$.
A probability measure $\mu\in\Cal P(X)$  is {\it
$\r$-non-singular} if $A\in\mathcal B,\ \mu(A)=0\ \Rightarrow\
\mu(\r(A))=0$ where $\r(A):=\{y\in X:\ \exists\ x\in A,\
(x,y)\in\r\}$. If $\r=\r_G$ as above where
$G$ is a countable group of measurable, invertible
transformations of $X$, then
$\r(A)=\bigcup_{g\in G}g(A)$ and we see that $\r(A)\in\mathcal B$. \par As shown in \cite{FM}, if a measure $\mu$ is
$\r$-non-singular then $\exists\ D=D_{\r,\mu}:\r\to\Bbb R_+$
measurable so that any holonomy $\Phi:B\overset{\r}\to\rightarrow
C$ is $\mu$-non-singular with
$\tfrac{d\mu\circ\Phi}{d\,\mu}(x)=D(x,\Phi(x))$ for $\mu$-a.e.
$x\in B$.

Let $\Psi:\r\to\Bbb R_+$ be measurable. A measure $\mu$ is called
$(\Psi,\r)$-{\it conformal} if it is $\r$-non-singular and
$D_{\r,\mu}\equiv \Psi$.  The  measure $\mu$ is called $\r$-{\it
invariant} if it is $(1,\r)$-{\it conformal} (i.e.
$D_{\r,\mu}\equiv 1$). We denote the collection of $(\Psi,\r)$-
conformal probabilities on $X$ by $\Cal P(X,\r,\Psi)$ and the
collection of $\r$-invariant probabilities on $X$ by $\Cal
P(X,\r)$.

If $T$ is an action of $\Bbb Z^d$ on $X$, and $\mu\in\Cal
P(X,T):=\Cal P(X,\r_T)$, then $(X,T,\mu)$ is called a {\it $\Bbb
Z^d$-random field}. \subheading{0.3 Ergodicity, transitivity,
irreducibility and mixing}

A measure $\mu$ is {\it $\r$-ergodic} if $\forall\ U,\ V\in\B,\
\mu(U),\mu(V)>0$, $(U\x V)\cap\r\ne\emptyset$; equivalently: the
collection of measurable, $\mathcal R$-invariant sets
$$\goth I(\r):=\{A\in\B:\ (x,y)\in\mathcal R\ \Rightarrow\ y\in A\}
\overset{\mu}\to{=}\{\emptyset,X\}.$$ We denote the collection of
$(\Psi,\r)$- conformal, $\r$-ergodic probabilities on $X$ by $\Cal
P_e(X,\r,\Psi)$ and the collection of $\r$-invariant, $\r$-ergodic
probabilities on $X$ by $\Cal P_e(X,\r)$.

The equivalence relation $\r$ on $X$ is {\it topologically
transitive} if $\forall\ U,\ V\subset X$ open, nonempty, $(U\x
V)\cap\r\ne\emptyset$.

\

Let  $X$ be a $\Bbb Z^d$-subshift. The shift action  $(X,T)$  is:\bul {\it topologically
transitive}
 if for any $U,V\subset X$ open, nonempty,
$\exists\ k\in\Bbb Z^d$ so that $U\cap T_k V\ne\emptyset$; and
\bul {\it topologically
mixing}
 if for any $U,V\subset X$ open, nonempty,
$\exists\ F\subset\Bbb Z^d$ finite, so that $U\cap T_k V\ne\emptyset\ \forall\  k\in\Bbb Z^d\setminus F$.

\

 A $\Bbb Z^d$-TMS $X$ is:\bul {\it irreducible}
 if the shift action $(X,T)$  is topologically
transitive; \bul {\it mixing} if if the shift action $(X,T)$  is topologically mixing; and
\bul
 {\it strongly irreducible} if  $\exists\ r>0$ so that
whenever $F,\ G\subset\Bbb Z^d$ finite, $\|j-k\|\ge r\ \forall\
j\in F,\ k\in G$:
$$X_{F\cup G}\cong X_F\x X_G.$$
\subheading{0.4 Tail conformality}
\par A {\it local potential} is a function  $G:X_{B(0,r)}\to\Bbb R$ where $X$ is a $\Bbb Z^d$-TMS and $r\ge 0$.
A local potential
 is called:\bul a {\it nearest neighbour -} or {\it Markov potential} in case $r=1$, and
 \bul   a {\it site potential} (or {\it activity function}) in case $r=0$ (i.e. $G:S\to\Bbb R$).

\

 Given a local potential $G:X_{B(0,r)}\to\Bbb R$, we call $p\in\mathcal P( X )$ \ \ {\it $G$
 - conformal }  if $p$ is
 $(e^{\Psi_G},\goth T( X ))$-conformal where
 $\Psi_G(x,y):=\sum_{j\in\Bbb Z^d}(G((T_jy)_{B(0,r)})-G((T_jx)_{B(0,r)}))$.

 The measure $p\in\mathcal P( X )$ is called
{\it locally conformal}, {\it nearest neighbour conformal} or {\it
site conformal } and if it is $G$ -conformal  for some local,
nearest neighbour or site potential $G$ (respectively).
\par Conformal measures with more general (e.g.
H\"older continuous) potentials are considered e.g. in \cite{R}, \cite{PS} and
\cite{S2}, where they are called {\it Gibbs measures}.
\par Evidently, any site conformal measure is exchangeable
($\mathcal E( X )$-invariant). Here, we show that global (i.e. globally
supported), ergodic exchangeable measures on certain kinds of $\Bbb
Z^d$-topological Markov shifts (see below) are forced to be site
conformal.
\par The  well known ``thermodynamic limit theorem" (see \cite{R} and also
\cite{PS},\cite{S2}) ensures the existence of  conformal measures for a $\Bbb Z^d$-TMS as weak limits of apporpriate sequences of atomic measures. \proclaim{Thermodynamic limit theorem} If
$X$ is a $\Bbb Z^d$-TMS and $G$ is a local potential, then
$\exists$ a $G$-conformal measure $p\in\mathcal P( X )$.
\endproclaim
\subheading{0.5 Equilibrium measures and tail conformal  measures}
 The {\it entropy} of a measure
$\mu\in\Cal P(X,T)$ (where $T$ is the shift $\Bbb Z^d$-action on
the  $\Bbb Z^d$-subshift $X$) is
$$h_\mu(X,T):=\lim_{n\to\infty}\tfrac1{|B(0,n)|}\sum_{x\in
X_{B(0,n)}}\mu([x]_{B(0,n)})\log\tfrac1{\mu([x]_{B(0,n)})}$$ (the
limit on the right hand side exists due to subadditivity). For any
$\mu\in\Cal P(X,T)$,
$$h_\mu(X,T)\le\lim_{n\to\infty}\frac{\log
|X_{B(0,n)}|}{|B(0,n)|}.$$ The limit on the right hand side (which
exists due to subadditivity) is called the {\it topological
entropy} of $(X,T)$ and is denoted $h(X,T)$.
\par Let $\mu\in\Cal P(X,T)$. The  $\mu$-{\it pressure} of a local potential $G$   on $X$ is
$P_\mu(G,X,T):=h_\mu(X,T)+\int_XGd\mu$.  The measure $\mu\in\Cal
P_e(X,T)\}$ is called a {\it $G$-equilibrium measure} (or {\it
equilibrium measure for $G$-pressure}) if $P_\mu(G,X,T)$ is
maximal. For existence of $G$-equilibrium measures in this situation,  (see \cite{R} and also \cite{M}).

Note that  $P_\mu(0,X,T)=h_\mu(X,T)$ and a $0$-equilibrium measure
is a measure of maximal entropy. We shall sometimes refer to a measure of maximal entropy as an
{\it equilibrium measure} (suppressing the $0$).
\par Let  $G$ be a  local potential on a TMS $X$.
\bul By the generalised Lanford-Ruelle theorem (the first part of theorem 4.2 in  \cite{R}, see also theorem 1.19 in \cite{BS2} and the original \cite{LR}),
 any $G$-equilibrium measure is $G$-conformal.
\bul By the generalised Dobrushin theorem (the second part of theorem 4.2 in  \cite{R}, see also proposition 4.1 in \cite{BS2} and the original \cite{D}),
if $X$ is strongly irreducible, then  any $G$-conformal
$\mu\in P_e(X,T)$ is a $G$-equilibrium measure.

 \subheading{0.6 Skew products and cocycle subrelations}
\par Let $\r$ be a countable equivalence relation on  $X$ and let
$\Bbb G$ a locally compact, Polish, Abelian (LCAP) topological
group. A Borel function $\Psi:\Cal R\to\Bbb G$ is called a $\Cal R
$-{\it cocycle} if
$$\Psi(x,z)=\Psi(x,y)+\Psi(y,z)\text{ whenever }(x,y),\
(y,z)\in\r.$$ For example, if $\mu\in\Cal P(X)$ is
$\r$-non-singular, then $\log D_{\mu,\r}:X\to\Bbb R$ is a
$\r$-cocycle.

 In this situation, we
consider the $\Psi$-{\it skew product relation}:
  $$
\r_\Psi:=\bigl\{\bigl((x,t),(y,s)\bigr)\in (X\x\Bbb G)^2: (x,y)\in
\r\text{ and }t-s=\Psi(x,y)\bigr\};
$$
 and $\Psi$-{\it subrelation}:$$ \r[\Psi]:=\{(x,y)\in X\x X: (x,y)\in\r\text{
and }\Psi(x,y)=0\}
$$ and note that
$$
\r_\Psi\cap (X\x\{0\})^2=\bigl\{\bigl((x,0),(y,0)\bigr)\in
(X\x\{0\})^2: (x,y)\in \r[\Psi]\bigr\}
$$
In case $X$ is a topological space and the countable equivalence
relation $\r$ is topological, we call the $\r$-cocycle $\Psi:\Cal
R\to\Bbb G$ {\it topological} if $\forall$ topological
$\r$-holonomy $\pi:A\overset{\r}\to\rightarrow B$,
 $x\mapsto\psi(x,\pi(x))$ is a continuous map $A\rightarrow\Bbb G$.
\par Let $X\subset
 S^{\Bbb Z^d}$ ($d\ge 1$) be a TMS, let $\Bbb G$ be a countable, Abelian group and let $G:S\to\Bbb G$.
  Define  $\Psi_G:\goth T(X)\to\Bbb G$  by
 $\Psi_G(x,y):=\sum_{j\in\Bbb Z^d}(G(y_j)-G(x_j))$. In the
 notation established above,
 $$\Cal E(X)=\goth T(X)[\Psi_{\sharp}]$$ where
 ${\sharp}:S\to\Bbb
G=\Bbb G_S:=\Bbb Z^{S}$ is defined by ${\sharp}(s):=e_s$ and
$(e_s)_t:=\d_{s,t}$.

\heading \S1 Main Results\endheading \subheading{1.1 Tail
non-singular exchangeable measures on a TMS}

\proclaim{Theorem 1} Suppose that $X\subset S^{\Bbb Z^d}$ is an
irreducible TMS with $\goth T(X)$ topologically transitive. Let
$G:G\to\Bbb G$, a countable Abelian group.

If $p\in\mathcal P(X)$ is global, $\goth T(X)$-nonsingular,  and $\goth
T(X)[G]$-invariant, ergodic, then $p$ is $H\circ G$-conformal for
some homomorphism $H:\Bbb G\to\Bbb R$.\endproclaim \subheading{Remarks}

1) Theorem 1 (and theorem 2 below) give information on exchangeable measures in case $G=\sharp:S\to\Bbb Z^S$ as in 0.6.

2)  Theorem 1 is a
partial converse to proposition 3.3 in \cite{S2}, which shows that
a $\goth T(X)$-non-singular, ergodic measure is $\goth
T(X)[G]$-ergodic.\subheading{1.2 Exchangeable measures on a
strongly aperiodic TMS} Let $X$ be a  $\Bbb Z^d$-TMS  and let
$\Bbb G$ be a countable Abelian group. We call a site function
$G:S\to\Bbb G$ {\it strongly aperiodic } (with respect to $X$) if
for every subgroup
$$\Bbb K\lvertneqq\Bbb H_G=\Bbb H_{X,G}:=\<\{\Psi_G(x,y):\ (x,y)\in\goth T(X)\}\>,$$
$$\exists\ F\subset\Bbb Z^d\ \ \text{\rm s.t.}\ \forall\  a\in X_F\
\exists\ b=b_a\in X_F,\ a_{\bdy F}=b_{\bdy F},\
\Psi_G(b,a)\notin\Bbb K$$ where $\<B\>$ denotes the group
generated by $B$. We call $X$ {\it strongly aperiodic} if every
site function is strongly aperiodic (with respect to $X$).

Examples of strongly aperiodic TMS's are given in \S3.

\proclaim{Theorem 2} Let  $X\subset S^{\Bbb Z^d}$ be a strongly
aperiodic, irreducible $\Bbb Z^d$-TMS, let $\Bbb G$ be a countable
Abelian group and let $G:S\to\Bbb G$. If $p\in\Cal P(X)$ is
global, $\goth T(X)[\Psi_G]$-invariant, ergodic, then $p$ is
$H\circ G$-conformal for some homomorphism $H:\Bbb G\to\Bbb
R$.\endproclaim \subheading{Remarks}

1) The existence of a global, $\mathcal E(X)$-invariant, ergodic
measure implies that $\goth T(X)$ topologically transitive.

2) If $X$ is strongly aperiodic, then any global, exchangeable
ergodic probability is site conformal.

3) As shown in proposition 3.1 (see \S3), any mixing $\Bbb Z$-TMS
is  strongly aperiodic.

4) Corollary 2.1 below is a stronger version of theorem 2 in case
$d=1$. This  generalizes theorem 6.2  of \cite{PS} (dispensing
with the assumption of shift-invariance). \subheading{1.3 Spectral
abundance of equilibrium measures}
\par Let $(\Om,\Cal F,P,T)$ be an ergodic, $\Bbb Z^d$-random field.
  A $T$-{\it
eigenfunction} is a function $f\in L^2$ satisfying $f\circ
T_k=e^{2\pi i\<\th,k\>}f\ \ (k\in\Bbb Z^d)$ for some $\th\in\Bbb
T^d$ (called the {\it eigenvalue}). The random field is called
\bul {\it totally ergodic} if each transformation $T_k\ (k\ne 0)$ is ergodic
(equivalently there are no $T$-eigenfunctions with rational
eigenvalues) and \bul{\it weakly mixing} if there are no non
trivial $T$-eigenfunctions.
\bul {\it mildly mixing} if there are no non
trivial {\it $T$-rigid sets}, a set $A\in\B$ being $T$-{\it rigid } if $\exists\ n_k\in\Bbb
Z^d,\ n_k\to\infty$ so that $P(T_{n_k}A\D A)\to 0$; and
\bul {\it (strongly) mixing} if  $P(A\cap T_nB)\to P(A)P(B)$ as $n\to\infty\ \forall\ A,B\in\Cal F$.

\par These (progressively stronger) properties are {\it spectral} properties in that they depend only on the
{\it spectral measure type } of  $(\Om,\Cal F,P,T)$: i.e. the measure class of $\sigma\in\Cal P(\Bbb T^d)$ defined by
the property:
$$\forall\ \nu\ll\s,\ \exists\ f\in L^2(P),\ \int_\Om ff\circ T_kdP=\widehat\nu(k).$$
See \cite{N}.

Now let $X$ be a $\Bbb Z^d$-subshift. We call a collection $\Cal
Q\subset\Cal P_e(X,T)$ {\it spectrally abundant} if  there are (different)
 measures $\mu\in\Cal Q$ so that:\bul  $(X,\mu,T)$ is not
totally ergodic;\bul  $(X,\mu,T)$ is totally ergodic but not weakly mixing;\bul
 $(X,\mu,T)$ is weakly mixing but not mildly mixing; \bul  $(X,\mu,T)$ is mildly mixing but
not strongly mixing;\bul  $(X,\mu,T)$  is strongly mixing.

 \proclaim{Theorem 3 }\par 1) There exists a
strongly aperiodic, strongly irreducible $\Bbb Z^3$-TMS whose
collection of
 equilibrium
measures is spectrally  abundant. \par 2) There is a nearest
neighbour potential $G$ on $\{0,1\}^{\Bbb Z^3}$ whose collection
of
 $G$-equilibrium
measures is spectrally  abundant.
\endproclaim

\heading \S2 Topological equivalence relations and the proof of
theorem 1
\endheading For $\G$ countable, $S$ finite$, X\in\B(S^\G)$, and $G:S\to\Bbb G$ (a countable, Abelian group) let
$$\Bbb H=\Bbb H_{X,G}:=\<\{\Psi_G(x,y):\ (x,y)\in\goth T(X)\}\>.$$
\proclaim{Skew product lemma} Suppose that $X\subset S^\G$ is
closed and that $\goth T(X)_{\Psi_G}$ is topologically transitive
on $X\x\Bbb H$.

If $p\in\mathcal P(X)$ is $\goth T(X)$-nonsingular, $\goth
T(X)[G]$-invariant, ergodic, then $p$ is site
conformal.\endproclaim\demo{Proof}

\par There is a unique $\s$-finite measure $m\in\goth M(X\x\Bbb H)$, $\goth T(X)_{\Psi_G}$-invariant,
ergodic such that $m(A\x\{0\})=p(A)$.

For $g\in\Bbb H$, let $q_g(A):=m(A\x\{g\})$. We claim first that
$q_g\ll p\ \forall\ g\in\Bbb H$.

To see this, let $A\in\B(X),\ q_g(A)>0$, then by $\goth
T(X)_{\Psi_G}$-ergodicity, $\exists\ A'\subset A,\ q_g(A')>0$ and
$\overline \pi:A'\x\{g\}\overset{\goth T(X)_{\Psi_G}}\to\to
\overline \pi(A'\x\{g\})=:\pi(A')\x\{0\}$. By $\goth
T(X)_{\Psi_G}$-invariance,
$$p(\pi(A'))=m(\pi(A')\x\{0\})=m(A'\x\{g\})=q_g(A')>0.$$ Evidently
$\pi:A'\overset{\goth T(X)}\to\to \pi(A')$ whence by $\goth
T(X)$-nonsingularity of $p$, $p(A')>0$. Thus $q_g\ll p\ \forall\
g\in\Bbb H$.

Next, set for $g\in\Bbb H,\ Q_g(x,y):=(x,y+g)\ \ \ (Q_g:X\x\Bbb
H\to X\x\Bbb H)$, then $m\circ Q_g$ is also $\goth
T(X)_{\Psi_G}$-invariant, ergodic, whence either $m\circ Q_g\sim
m$ or $m\circ Q_g\perp m$.

Let $\Bbb K:=\{g\in\Bbb H:\ m\circ Q_g\sim m\}$, then $\Bbb K$ is
a subgroup of $\Bbb H$ and $q_g\sim p$ if $g\in\Bbb K$ and
$q_g\equiv\ 0\ \ (q_g\ll p\ \&\  q_g\perp p)$ if $g\notin\Bbb K.$

 To see that $\Bbb K\supseteq\Bbb H$ fix $g\in\Bbb H$. By
topological transitivity of $\goth T(X)_{\Psi_G}$, $\exists\
A\in\B(X),\ p(A)>0$ and $\pi:A\x\{0\}\overset{\goth
T(X)_{\Psi_G}}\to\to \pi(A\x\{0\})\subset X\x\{g\}$. Thus
$$q_g(X)\ge m(\pi(A\x\{0\}))=m(A\x\{0\})=p(A)>0,$$
$q_g\neq 0$ and $g\in\Bbb K$.

It now follows that $\exists$ a homomorphism $ H:\Bbb H\to\Bbb R$
so that $m\circ Q_g=e^{H(g)}m$ whence
$$\tfrac{dp\circ\pi}{dp}(x)=e^{H(\Psi_G(x,\pi(x)))}\ \
\forall\ \pi:A\overset{\goth T(X)}\to\to \pi(A)$$ and $p$ is site
conformal.

 \hfill\qed\enddemo

Suppose that $\psi:\r\to\Bbb G$ is a topological $\r$-cocycle
($\r$ a topologically transitive, topological equivalence
relation, $\Bbb G$ a LCAP topological  group). We call $g\in \Bbb
G$  a {\it topological essential value}  of $\psi$ if
$$(\psi^{-1}U_g)\cap (A\x A)\ne\emptyset\ \ \forall
\  A\subset X,\ U_g\subset\Bbb G\ \ \text{\rm nonempty,\ open sets
with }\ g\in U_g.$$ The collection of topological essential values
of $\psi$ is denoted
$$\Bbb E=\Bbb E_{\text{\rm top}}(\psi)$$
and forms  a closed subgroup of $\Bbb G$ (see \cite{LM}).
\par We need the following  version of proposition 3.2 in
\cite{LM}. \proclaim{Topological essential value lemma} Suppose
that $\r$ is topologically transitive on $X$, then $\r_{\psi}$ is
topologically transitive on $X\x\Bbb G$ $\iff$ $\Bbb E_{\text{\rm
top}}(\psi)=\Bbb G$.\endproclaim \demo{Proof}\f$\Rightarrow$)\ \
Suppose that $\r_{\psi}$ is topologically transitive on $X\x\Bbb
G$, $g\in\Bbb G$ and $A\subset X,\ g\in U_g\subset\Bbb G$ are
nonempty and open. There is an open neighborhood $V$ of $0$ so
that $g+V-V\subset U_g$. By definition $((A\x V)\x
(A\x(g+V))\cap\r_\psi\ne\emptyset$, whence
$$(\psi^{-1}U_g)\cap (A\x A)\supset(\psi^{-1}(g+V-V))\cap (A\x A) \ne\emptyset.$$ Thus $g\in\Bbb
E$. \smallskip\f$\Leftarrow$)\ \   Now suppose that $\r_\psi$ is
not topologically transitive on $X\x\Bbb G$, then $\exists\ A,\
B\subset X$, $U,\ V\subset\Bbb G$ open with $0\in U\cap V$ and
$g\in \Bbb G$ so that $\r_\psi\cap ((A\x U)\x(B\x
(g+V)))=\emptyset$.
\par By topological transitivity of $\r$ on $X$, $\exists\ A'\subset
A,\ B'\subset B$ open and a topological holonomy
$\pi:A'\overset{\r}\to\rightarrow B'$ so that
$x\mapsto\psi(x,\pi(x))$ is continuous.
\par Fix $0\in W\subset\Bbb G$ so that
$W+W\subset V$. Using continuity of $x\mapsto\psi(x,\pi(x))$ we
ensure (by possibly reducing $A',\ B'$) that $\exists\ h\in\Bbb G$
with $\psi(x,\pi(x))\in h+W\ \forall\ x\in A'$.
\par We claim that $k:=g-h\notin\Bbb E$. To see this, note first
that
$$\r_\psi\cap ((A'\x U)\x(A'\x
(k+W)))=\emptyset.$$ Otherwise $\exists\ (x,y)\in A'\x U,\
(x',y')\in A'\x (k+W)$ with $((x,y),(x',y'))\in\r_\psi$, whence
$((x,y),(\pi(x'),y'+\psi(x',\pi(x'))\in\r_\psi$. However
$$\align (\pi(x'),y'+\psi(x',\pi(x')) & \in B'\x
(k+W+\psi(x',\pi(x'))\\ &\subset B'\x (k+W+h+W)\\ &\subset B'\x
(g+V)\endalign$$ contradicting $\r_\psi\cap ((A\x U)\x(B\x
(g+V)))=\emptyset$.
\par To finish the proof that $k\notin\Bbb E$, fix $0\in W_0\subset U$
open so that $W_0+W_0\subset W$. If $k\notin\Bbb E$, then
$\exists\ (x,x')\in\r\cap (A'\x A')$ with $\psi(x,x')\in k+W_0$.
It follows that for $y\in W_0,\ (x',y+\psi(x,x')\in A'\x
(k+W_0+W_0)\subset (k+W)$ contradicting $\r_\psi\cap ((A'\x
U)\x(A'\x (k+W)))=\emptyset.$
 \hfill\qed\enddemo
\proclaim{Transitivity lemma} Let $X\subset
 S^{\Bbb Z^d}$ be an irreducible TMS such that $\goth T(X)$ is  topologically transitive on $X$.
 \par Suppose that $G:S\to\Bbb G$, then
$$\Bbb E_{\text{\rm
top}}(\Psi_G)=\overline{\<\{\Psi_G(x,y):\ (x,y)\in \goth
T(X)\}\>}.$$
 \endproclaim \demo{Proof of $\supseteqq$ (as in \cite{S2})}
Let $(x,y)\in\goth T(X),\ \Psi_G(x,y)=g$. We show that $g\in\Bbb
E$. To this end,  let $F\subset\Bbb Z^d$ be finite, $a\in X_F$.
We'll show that $\exists\ \Lambda\subset\Bbb Z^d$ finite,
$\Lambda\supset F$ and $b,c\in X_\Lambda,\ b_F=c_F=a,\
b_{\bdy\Lambda}=c_{\bdy\Lambda}$ so that $\Psi_G(u,v)=g\ \forall\
u\in [b],\ v\in [c],\ u_{\Lambda^c}=v_{\Lambda^c}$.

\par Since $(x,y)\in\goth T(X),\ \exists\ B\subset\Bbb Z^d$  a
cube so that $x_{B^c}=y_{B^c}$. WLOG, $x_{\bdy B}=y_{\bdy B}$. By
irreducibility, $\exists\ k\in\Bbb Z^d$ so that $F\cap
(B+k)=\emptyset$ and $[a]_F\cap T_k[x_B]_B\ne\emptyset.$ Now let
$\Lambda:=F\cup (B+k),$ then $\bdy\Lambda=\bdy F\cup(k+\bdy B)$.

Define $b\in X_\Lambda$ by $b=z_\Lambda$ where $z\in [a]_F\cap
T_k[x_B]_B.$ Evidently $b_F=a$ and $b_{j}=x_{j-k}\ \forall\ j\in
B+k$. Now define $c\in S^\Lambda$ by
$$c_j=\cases y_{j-k}\ \ \ \ \ j\in B+k,\\ b_j\ \ \ \ \text{\rm
else.}\endcases$$ Since $x_{\bdy B}=y_{\bdy B}$ we have that $c\in
X_\Lambda$ and $b_{\bdy\Lambda}=c_{\bdy\Lambda}$.

It follows that $\forall\ u\in [b],\ v\in [c],\
u_{\Lambda^c}=v_{\Lambda^c}$,
$$\Psi_G(u,v)=\Psi_G(x,y)=g.$$
\hfill\qed\enddemo \demo{Proof of theorem 1}\ \  By the
transitivity lemma and topological essential value lemma, $\goth
T(X)_{\Psi_G}$ is topologically transitive on $X\x\Bbb H$. Theorem
1 now follows from the skew product lemma.\hfill\qed\enddemo

\par Let $X$ be a topologically transitive $\Bbb Z$-TMS, then  (see e.g. \cite{Ch})

$X=\biguplus_{k=0}^{N-1}X_k$ where $N\in\Bbb N$ and
$X_0,\dots,X_{N-1}$ are disjoint, clopen subsets of  $X$ with
$TX_k=X_{k+1\ \mod N}$; and each $(X_k,T^N)$ is mixing.

 This
decomposition is called the {\it periodic decomposition} of $X$ (also known as the {\it cyclic} or {\it spectral decomposition}),
$N=N^{(X)}$ is called the {\it period} of $X$ and each $X_k$ is
called a {\it basic, mixing} set for $X$.
\par Note that each $X_k$ is $\goth T(X)$-invariant. By theorem 3.3 of \cite{PS},
any globally supported Markov measure on $X_k$ is $\goth
T(X_k)[\Psi_G]$-nonsingular, ergodic. In particular, if $H:\Bbb
G\to\Bbb R$ is a homomorphism and $\mu$ is the  $(X_k,T^N)$-Gibbs
measure on $X_k$ with potential $e^{H\circ G_N}$ where $G_N:=\sum_{k=0}^{N-1}G\circ T^k$ (unique,
$T^N$-invariant), then $\mu$ is $\goth T(X)[\Psi_G]$-invariant,
ergodic.
 \proclaim{Corollary 2.1}

Let  $X\subset S^{\Bbb Z}$ be a topologically transitive $\Bbb
Z$-TMS, let $\Bbb G$ be a countable Abelian group and let
$G:S\to\Bbb G$.

If $p\in\Cal P(X)$ is
 $\goth T(X)[\Psi_G]$-invariant, ergodic, then

 then
there exist \roster\item a homomorphism $H:\Bbb G\to\Bbb R$;\item
a $\goth T(X)[\Psi_G]$-invariant, topologically transitive,  TMS
$X'\subset X$;
\item a basic mixing set  $X_0'\subseteq X'$;
\endroster
so that $\text{\rm supp}\, p=X_0'$  is the $T^N$-Gibbs measure on
$X_0'$ with potential $e^{H\circ G_N}$. \endproclaim \demo{Proof} We
claim first that for each $s\in S$, either $N_s(x):=\sum_{n\in\Bbb
Z}\d_{x_n,s}=0\ \ p$-a.s., or $N_s(x)=\infty \ p$-a.s.. This is
because if $k\in\Bbb N,\ p([N_s=k])>0$ then $\exists\
K_n\subset\Bbb Z,\ |K_n|=k\ \ (n\ge 1)$ so that the sets
$$A_n:=\{x\in X:\ x_j=s\ \forall\ j\in K_n,\ x_k\ne s\ \forall\ j\notin K_n\}$$
are disjoint, exchangeably- (whence $\goth T(X)[\Psi_G]$-)
equivalent and thus with equal positive measure, entailing
$p(X)=\infty$.
\par Next, as in step 1 of the proof of  theorem 5.0 of \cite{ANS}, $p$ is the restriction
of an irreducible, shift invariant, Markov measure $\mu\in\Cal
P(X)$ to some clopen set in $X$.
\par Let $X'=\text{\rm supp}\, \mu$, then $X'$ is a topologically transitive TMS.

Let $X'=\biguplus_{k=0}^{N-1}X_k'$ be the periodic decomposition
 of $X'$, by $\goth T(X)[\Psi_G]$-ergodicity of $p,\ \exists\ k$ so that
 $\text{\rm supp}\, p\subseteq X_k'$. Since $\mu$ is $\goth T(X_k')[\Psi_G]$-nonsingular, ergodic,
 $\text{\rm supp}\, p= X_k'$. The result now follows from theorem 1.
\hfill\qed\enddemo
\heading{\S3 Conditions for strong aperiodicity and the proof of
theorem 2 }\endheading\proclaim{Proposition 3.1}
\par a) A $\Bbb Z^d$-TMS $X\subset S^{\Bbb Z^d}$ is strongly
aperiodic iff\ \ $\sharp:S\to\Bbb Z^S$ is strongly aperiodic.\par
b)  Any mixing $\Bbb Z^1$-TMS is strongly
aperiodic.\endproclaim\demo{Proof}\par a) Suppose that $X$ is
strongly aperiodic, let $\Bbb G$ be a countable Abelian group and
let $G:S\to\Bbb G$ be a site function.

Define $\pi=\pi_G:\Bbb Z^S\to\Bbb G$ by $\pi(\sum_{s\in
S}n_se_s):=\sum_{s\in S}n_sG(s)$ then $\pi$ is a homomorphism and
$\pi\circ\sharp=G$ whence $\Bbb H_G=\pi_G\Bbb H_\sharp$. Strong
aperiodicity of $G$ follows from this.

\par b)  Let $G:S\to\Bbb G$ be a site function and consider
$g:X_+:=\{x_+=(x_1,x_2,\dots):\ x=(\dots,x_{-1},x_0,x_1,\dots)\in
X\}\to\Bbb G$ be defined by $g(x):=G(x_1)$. By the well-known
cohomology lemma (see e.g. lemma 4.3 in \cite{ANS}), $g=a+h-h\circ
T+\overline g$ where
  $a\in\Bbb G,\ \overline g:X_+\to \Bbb G_g:=\<\{g_n(x)-g_n(x'):\ n\ge 1,\ T^nx=x,\ T^nx'=x'\}\>$  and
$h:X\to\Bbb G$ are both generated by site functions such that
$\overline g :X\to \Bbb G_g$ is aperiodic (in the sense that $\Bbb
G_{\overline g}=\Bbb G_g$). It follows that $\Bbb H_G=\Bbb G_g$
and that $\Psi_{\overline g}\equiv\Psi_G$.
\par As mentioned in the proof of theorem 2.2 in \cite{ANSS},
$\forall\ H\lvertneqq \Bbb G_g,\ \exists\ \ell\ge 1$ so that
$\forall$ $a=(a_1,\dots,a_{\ell+1})\in X_{\{1,\dots,\ell+1\}},\
\exists$ a path $b=b_a=(b_1,\dots,b_{\ell+1})\in
X_{\{1,\dots,\ell+1\}}$ such that $a_1=b_1,\
a_{\ell+1}=b_{\ell+1}$ and $f_\ell(a)-f_\ell(b)\notin H$. This is
strong aperiodicity of $G$.\hfill\qed\enddemo \demo{Proof of
theorem 2}

There is a unique $\s$-finite, $\goth T(X)_{\Psi_{ G
}}$-invariant, ergodic measure $m\in\goth M(X\x\Bbb H)$ such that
$m(A\x\{0\})=p(A)$ (where $\Bbb H=\Bbb H_G$ as above).

Set for $g\in\Bbb H,\ Q_g(x,y):=(x,y+g)\ \ \ (Q_g:X\x\Bbb H\to
X\x\Bbb H)$, then $m\circ Q_g$ is also $\goth T(X)_{\Psi_{ G
}}$-invariant, ergodic, whence either $m\circ Q_g\sim m$ or
$m\circ Q_g\perp m$.

\f$\bullet\ \ \  $ Let $\Bbb K:=\{g\in\Bbb H:\ m\circ Q_g\sim
m\}$, then $\Bbb K$ is a subgroup and $q_g\sim p$ if $g\in\Bbb K$.
We'll show that $\Bbb K=\Bbb H$. Assume otherwise.

\f$\bullet\ \ \ $ Let $F\subset\G$ be as the definition of strong
aperiodicity  adapted to $\Bbb K$ and let
$$J:=\{\Psi_ G (b_a,a):\ a\in X_F\}\subset\Bbb H\setminus\Bbb K.$$ For $g\in\Bbb H$, let $q_g(A):=m(A\x\{g\})$ and set
$\overline q:=\sum_{j\in J}q_j$, then $\overline  q\perp p$.

 \f$\bullet\ \ \ $ {\bf Claim }  $\exists\ E\subset\G,\ \goth a\in X_E$ so that $\overline  q([\goth a])<\infty$.
\smallskip\f{\it Proof of claim}\ \ \  By irreducibility, $\exists\
\{k(a):\ a\in X_F\}\subset\G$ so that $\{F+k(a)\}_{a\in X_F}$ are
pairwise disjoint and $\exists\  \goth a\in X_E\ \
(E=\biguplus_{a\in X_F}(F+k(a)))$ so that
$$\goth a_i=a_{i-k(a)}\ \forall\ i\in F+k(a),\ \ a\in X_F.$$
It suffices to prove that $q_j([\goth a])\le 1\ \forall\ j\in J$.
Suppose that $j=\Psi_{ G }(b,a)$ where $a,b\in X_F,\ a_{\bdy
F}=b_{\bdy F}$. Define $\goth a'\in X_E$ by
$$\ \ \ \goth a_i'\ \ \ =\ \  \cases & b_i\ \ \ \ \ \ \ i\in F+k(a),\\ &
\goth a_i\ \ \ \ \ \ \ \text{else,}\endcases$$ then $c'_{\bdy
E}=c_{\bdy E}$ and $\Psi_{ G }(\goth a,\goth a')=\Psi_{ G
}(a,b)=-j$. It follows that

$[\goth a]\x\{j\}\overset{\goth T(X)_{\Psi_{ G }}}\to\to [\goth
a']\x\{0\}$, whence
$$q_j([\goth a])=m([\goth a]\x\{j\})=m([\goth a']\x\{0\})=p([\goth a'])\le 1.\qed$$

Next,  \f$\bullet\ \ \ $ Since $p$ is global and $\overline q\perp
p$, $\exists\ K\subset [\goth a]$ compact so that $p(K)>0,\
\overline  q(K)=0$. \f$\bullet\ \ \ $ $\exists\ U$ open $K\subset
U\subset [\goth a]$  so that $\overline
q(U)<\tfrac{p(K)}{2|X_F|}$. \f$\bullet\ \ \ $ $\exists\
E'\subset\G,\ c\in X_{E'}$ so that $\overline
q([c]_{E'})<\tfrac{p([c]_{E'})}{2|X_F|}$. \f$\bullet\ \ \ $ Fix
$k\in\G$ with $E'\cap (F+k)=\emptyset$ and set $\Lambda:=E'\uplus
(F+k)$. For $a\in X_F,\ b\in X_{E'} $, let $\<a,b\>\in S^\Lambda$
be defined by
$$\ \ \ \<a,b\>_j\ \ \ =\ \  \cases & b_j\ \ \ \ \ \ \ j\in E',\\ &
a_{j-k}\ \ \ \ \ \ \ j\in F+k.\endcases$$ Not all $\<a,b\>\in
X_\Lambda$, however if $Y_b:=\{a\in X_F:\ \<a,b\>\in X_\Lambda\}$,
then $[b]=\biguplus_{a\in Y_b}[\<a,b\>]$. Thus:
 \f$\bullet\ \ \ $ $\exists\ a\in X_F$
such that
$$p([\<a,c\>])\ge \tfrac{p([c])}{|X_F|}.$$
\f$\bullet\ \ \ $ Now $[\<a,c\>]\x\{0\}\overset{\goth T(X)_{\Psi_{
G }}}\to\to [\<b_a,c\>]\x\{j\}$ where $j=\Psi_{ G }(b_a,a)\in J$.
Thus
$$\tfrac{p([c])}{|X_F|}\le
m([\<a,c\>]\x\{0\})=m([\<b_a,c\>]\x\{j\})\le \overline  q([c])<
\tfrac{p([c])}{2|X_F|}.$$ Thus $\Bbb K=\Bbb H$.

\f$\bullet\ \ \ $ It now follows that $\exists$ a homomorphism $
H:\Bbb H\to\Bbb R$ so that $m\circ Q_g=e^{H(g)}m$ whence
$$\tfrac{dp\circ\pi}{dp}(x)=e^{H(\Psi_{ G }(x,\pi(x)))}\ \
\forall\ \pi:A\overset{\goth T(X)}\to\to \pi(A)$$ and $p$ is site
conformal.

 \hfill\qed\enddemo

 \subheading{Condition
$\mho$} This condition implies strong aperiodicity (and is
 equivalent to it when $d=1$):
$$\exists\ F\subset\G\ \ \text{\rm s.t.}\
\forall\  a\in X_F,\ \<\{\Psi_\sharp(a,b):\  b\in X_F,\ a_{\bdy
F}=b_{\bdy F}\}\>=\Bbb H_{X,\sharp}.\tag{$\mho$}$$
\subheading{Condition $\maltese$}

TMS $X\subset S^{\Bbb Z^d}$ is such that $S=W\uplus Z$ with
$Z\ne\emptyset$ and so that \f$\bullet\ \ \ $ $\maltese$(i)
$\forall\ x\in X\ \ z\in Z, x'\in X$ where $x'_i=x_i\ \forall\
i\ne 0$ and $x_0'=z$; and \f$\bullet\ \ \
$$\maltese$(ii) If $a\in X_{[-1,1]^d}$ and $a_i\in Z\ \forall\
i\in [-1,1]^d\setminus\{0\}$ then $a^{(s)}\in X_F\ \forall\ s\in
S$ where $a^{(s)}_i=a_i\ \forall\ i\in [-1,1]^d\setminus\{0\}$ and
$a^{(s)}_0=s$. \subheading{Example: $\Bbb Z^2$ Iceberg model\ \
(\cite{BS2})} Here
$$S:=\{0,\pm 1,\pm 2,\dots,\pm M\},\ X:=\{x\in S^{\Bbb Z^2}:\ x_{n+e_i}x_n\ge 0\
\forall\ n\in\Bbb Z^2,\ i=1,2\}.$$ It is easy to see that the
iceberg model satisfies condition $\maltese$ with $Z=\{0\}$.

 \proclaim{ Proposition 3.2} $\maltese\ \ \Rightarrow$
$\mho$.\endproclaim\demo{Proof} Note that $\Bbb
H_{X,\sharp}\subseteq\Bbb H_{S^{\Bbb Z^d},\sharp}=\<\{e_s-e_t:\
s,\ t\in S\}\> \cong\Bbb Z^{S\setminus\{s_0\}}$ where $s_0\in Z$
is fixed.
  Define $\g:S\to\Bbb
Z^{S\setminus\{s_0\}}$ by $\g(s):=e_s$ for $s\ne s_0$ and
$\g(s_0)=0$. Choose $F:=[-2,2]^d$, then $\bdy F=F\setminus
F^\circ,\ F^\circ :=[-1,1]^d$.

Fix $a\in X_F$ and set $G(a):= \{\Psi_\g(a,b):\ b\in X_F,\ a_{\bdy
F}=b_{\bdy F}\}$. Set $g(a):=\sum_{i\in F^\circ}a_i$, then
$$G(a)\supset\{-\g(a_i):\ i\in F^\circ\}\cup\{-g(a)+e_s:\ s\in
S\}$$ whence  $\<G(a)\>=\Bbb Z^{S\setminus\{s_0\}}$ establishing
condition $\mho$. \hfill\qed\enddemo\subheading{Example:
generalized $\Bbb Z^d$ Beach
 model\ \ \ (\cite{BS2},\ \cite{H})} Here $S:=A\x B$ where $A=A_0\uplus A_1,\ A_0\ne\emptyset.$
 Writing $s\in S$ as
 $s=(\a (s),\b (s))\in A\x B$:
 $$X:=\{x\in S^{\Bbb Z^d}:\ \forall\ n\in\Bbb Z^d,\ 1\le i\le d,\
 \a (x_n), \a (x_{n+e_i})\in A_0\ \ \text{\rm or}\ \b (x_n)=\b (x_{n+e_i})
\}.$$
 \proclaim{ Proposition 3.3}
The generalized beach model satisfies
$\mho$.\endproclaim\demo{Proof} Recall that $S:=A\x B$ where
$A=A_0\uplus A_1,\ A_0\ne\emptyset$ and

 $$X:=\{x\in S^{\Bbb Z^d}:\ \forall\ n\in\Bbb Z^d,\ 1\le i\le d,\
 \a (x_n), \a (x_{n+e_i})\in A_0\ \ \text{\rm or}\ \b (x_n)=\b (x_{n+e_i})
\}$$ where $s=(\a (s),\b (s))\in A\x B=S.$ \par We note the
following facts: \bul Suppose that $x\in X,\ k\in\Bbb Z^d,\ s\in
S,\ \a (s)\in A_0,\ \b (s)=\b (x_k)$. If $y\in S^{\Bbb Z^d}$ is
defined by $y_{\Bbb Z^d\setminus\{k\}}=x_{\Bbb
Z^d\setminus\{k\}},\ y_k=s$ then $y\in X$; \bul Suppose that $a\in
X_{B_1(k,r)}$  (where $B_1(k,r):=\{k'\in\Bbb Z^d:\ \|k'-k\|_1\le
r\}$), then either $\exists\ j\in B_1(k,r)$ with $\a(a_j)\in A_0$
or $\exists\ b\in B$ so that $\b(a_j)=b\ \forall\ j\in B_1(k,r)$.

Now  let $F=B_1(0,6)$ and  let $a\in X_F$. We show that
$$\G_a:=\<\{\Psi_\sharp(a,b):\ b\in X_{B_1(j_0,3)},\ a_{\bdy
{B_1(j_0,3)}}=b_{\bdy {B_1(j_0,3)}}\}\>=\<\{e_z-e_w:\ z,\ w\in
S\}\>.$$

\f{\bf Case 1}\ \ If $\exists\ j_0\in B_1(0,3),\ \a(a_{j_0})\in
A_0$, we change $a_{B_1(j_0,3)}$ only.

Fix $z\in S$ and let $t_i\in S\ \ \ (i\in\bdy B_1(j_0,1))$ with
$\a(t_i)\in A_0,\ \b(t_i)=\b(z)$. Choose $s_j\in S,\ \a(s_j)\in
A_0,\ \b (s_j)=\b (a_j)\ \forall\ j\in\bdy B_1(j_0,2)$ and define
$b\in S^F$ by $$b_j=\cases & s_j\ \ \ \  \ \
\ j\in \bdy B_1(j_0,2),\\ & t_j\ \ \ \ \ \ \ j\in\bdy B_1(j_0,1),\\
& z \ \ \ \ \ \ \ j=j_0,\\ & a_j\ \ \ \ \ \ \
\text{else}.\endcases$$ It can be checked that $b\in X_F,\ a_{\bdy
{B_1(j_0,3)}}=b_{\bdy {B_1(j_0,3)}}$ and that
$$\Psi_\sharp(a,b)=\sum_{j\in\bdy
B_1(j_0,2)}(e_{s_j}-e_{a_j})+\sum_{i\in\bdy
B_1(j_0,1)}(e_{t_i}-e_{a_i})+e_z-e_{a_{j_0}}.$$

Other values of $\Psi_\sharp(a,b):\ b\in X_{B_1(j_0,3)},\ a_{\bdy
{B_1(j_0,3)}}=b_{\bdy {B_1(j_0,3)}}$ are obtained as follows: \bul
$\Psi_\sharp(a,b)=e_s-e_{a_j},\ j\in B_1(j_0,2),\ s\in S,\
\a(s)\in A_0,\ \b (s)=\b (x_j)$  where $b_j=s,\ b_i=a_i\ \forall\
i\ne j$; \bul $\Psi_\sharp(a,b)=\sum_{j\in\bdy
B_1(j_0,2)}(e_{s_j}-e_{a_j})+e_t-e_{a_i},\ \ s_j\in S,\ \a(s_j)\in
A_0,\ \b (s_j)=\b (a_j)\ \forall\ j\in\bdy B_1(0,2)$ and $i\in\bdy
B_1(0,1),\ t\in S,\ \a(t)\in A_0$ where $b_j=s_j\ \forall\ j\in
\bdy B_1(j_0,2),\ b_i=t,\ b_\ell=a_\ell\ \forall\ \ell\notin\bdy
B_1(j_0,2)\cup\{i\}$.
\par Thus we see that
$e_z-e_{a_{j_0}}\in\G_a\ \forall\ z\in S$ whence
$e_z-e_{z'}\in\G_a\ \forall\ z,z'\in S$. \f{\bf Case 2}\ \ If
$\nexists\ j_0\in B_1(0,3),\ \a(a_{j_0})\in A_0$, then $\exists\
\b\in B$ so that $\b(a_j)=\b\ \forall\ j\in B_1(0,3)$ and we
change $a_{B_1(0,3)}$ only.

Fix $z\in S$ and let $t_i\in S\ \ \ (i\in\bdy B_1(0,1))$ with
$\a(t_i)\in A_0,\ \b(t_i)=\b(z)$. Choose $s_j\in S,\ \a(s_j)\in
A_0,\ \b (s_j)=\b\ \forall\ j\in\bdy B_1(0,2)$ and define $b\in
S^F$ by $$b_j=\cases & s_j\ \ \ \  \ \
\ j\in \bdy B_1(0,2),\\ & t_j\ \ \ \ \ \ \ j\in\bdy B_1(0,1),\\
& z \ \ \ \ \ \ \ j=0,\\ & a_j\ \ \ \ \ \ \
\text{else}.\endcases$$ It can be checked that $b\in X_F,\ a_{\bdy
{B_1(0,3)}}=b_{\bdy {B_1(0,3)}}$ and that
$$\Psi_\sharp(a,b)=\sum_{j\in\bdy
B_1(0,2)}(e_{s_j}-e_{a_j})+\sum_{i\in\bdy
B_1(0,1)}(e_{t_i}-e_{a_i})+e_z-e_{a_{0}}.$$ Other values of
$\Psi_\sharp(a,b):\ b\in X_{B_1(0,3)},\ a_{\bdy
{B_1(0,3)}}=b_{\bdy {B_1(0,3)}}$ are obtained as follows: \bul
$\Psi_\sharp(a,b)=e_s-e_{a_j},\ j\in B_1(0,2),\ s\in S,\ \a(s)\in
A_0,\ \b (s)=\b$  where $b_j=s,\ b_i=a_i\ \forall\ i\ne j$; \bul
$\Psi_\sharp(a,b)=\sum_{j\in\bdy
B_1(0,2)}(e_{s_j}-e_{a_j})+\sum_{i\in B_1(0,1)}(e_{t_i}-e_{a_i}),\
\ s_j\in S,\ \a(s_j)\in A_0,\ \b (s_j)=\b\ \forall\ j\in\bdy
B_1(0,2)$ and $t_i\in S,\ \a(t_i)\in A_0\ \forall\ i\in B_1(0,1)$
where $b_j=s_j\ \forall\ j\in \bdy B_1(0,2),\ b_i=t_i\ \forall\
i\in B_1(0,1),\ b_\ell=a_\ell\ \forall\ \ell\in\bdy B_1(0,3)$.
\par Thus we see that
$\sum_{i\in B_1(0,1)}(e_{t_i}-e_{a_i})\in\G_a\ \forall\ t_i\in S,\
\a(t_i)\in A_0\ \ (i\in B_1(0,1))$ whence $e_z-e_{a_{j_0}}\in\G_a\
\forall\ z\in S$ and $e_z-e_{z'}\in\G_a\ \forall\ z,z'\in S$.
\hfill\qed\enddemo

\heading{\S4 Shift action on conformal measures, equilibrium measures and the proof of theorem 3
}\endheading

Suppose that  $X$ is a strongly irreducible $\Bbb Z^d$-TMS,  $G:X\to\Bbb R$ is a local
potential and  $\mu\in\Cal P(X,\goth T(X),\Psi_G)$. By the
ergodic decomposition (see \cite{F}, \cite{GS}),
$$\mu=\mu_\nu:=\int_{\Cal P_e(X,\goth T(X),\Psi_G)}\om d\nu(\om)$$
where $\nu\in\Cal P(\Cal P_e(X,\goth T(X),\Psi_G))$ and the
measure spaces \f$(\Cal P_e(X,\goth T(X),\Psi_G),\B(\Cal
P_e(X,\goth T(X),\Psi_G)),\nu)$ and $(X,\goth I(\goth T(X)),\mu)$
are isomorphic.

If $\om\in\Cal P_e(X,\goth T(X),\Psi_G)$ and $k\in\Bbb Z^d$, then
$\om\circ T_k \in\Cal P_e(X,\goth T(X),\Psi_G)$ and a Borel $\Bbb
Z^d$-action on $\Cal P_e(X,\goth T(X),\Psi_G)$ is defined by
$S_k\om:=\om\circ T_k$.

 \proclaim{Proposition
4.1}\par a) The measure $\mu\in\Cal P(X)$ is a $G$-equilibrium
measure iff $\mu=\mu_\nu$ where
$$\nu\in\Cal P_e(\Cal P_e(X,\goth T(X),\Psi_G),S).$$
\par In this case, \par b) any $T$-rigid set is $\goth
T(X)$-invariant; and \par c) $T$ is totally ergodic, weakly
mixing, mildly mixing iff $S$ has the respective property.

\endproclaim\demo{Proof} By theorem 4.2 in \cite{R},
the collection of $G$-equilibrium measures is given by $\Cal P_e(X,T)\cap \Cal P(X,\goth T(X),\Psi_G)$.

  Suppose first that $\mu\in\Cal P(X,T)$ is $\goth T(X)$-non-singular.
We claim that
$$\|f\circ
T_n\circ\pi\circ T_{-n}-f\|_1\to 0\ \ \forall\ f\in L^\infty,\
\pi\in [\goth T(X)].\tag{$\star$}$$ To see this, note first that
if $(x,y)\in\goth T(X)$ and $x_{F^c}=y_{F^c}$ where $F\subset \Bbb
Z^d$ is finite, then $(T_ny)_k=(T_nx)_k$ iff $n+k\notin F$. Thus
$\rho(T_nx,T_ny)\to 0$ as $n\to\infty$ where $\rho$ is a metric on
$X$ generating the standard (product) topology.

Now let $\pi\in [\goth T(X)]$. It follows that for $f\in C(X)$,
$|f(T_n\pi(x))-f(T_nx)|\to 0\ \forall\ x\in X$, whence $$\|f\circ
T_n\circ\pi-f\circ T_{n}\|_1\to 0.$$ Next, we obtain this
convergence for $f\in L^\infty$ by approximation. Let
$\phi:=\tfrac{d\mu\circ\pi^{-1}}{d\mu}$. Since $\mu\in\Cal
P(X,T)$, the sequence $\{\phi\circ T_n:\ n\in\Bbb Z^d\}$ is
uniformly absolutely continuous with respect to $\mu$.

Let $g\in L^\infty,\ |g|\le 1$ and let $\e>0$, then \bul $\exists\
\d>0$ so that if $A\in\B,\ \mu(A)<\d$, then $\sup_{n\in\Bbb
Z^d}\int_A\phi\circ T_nd\mu<\e$; and \bul $\exists\ f\in C(X),\
|f|\le 2$ so that $\mu([f\ne g])<\d$.

We see that as $n\to\infty$,
$$\align  \|g\circ T_n\circ\pi-g\circ
T_{n}\|_1 & \le
\|(g-f)\circ T_n\circ\pi\|_1+ \|(f-g)\circ T_{n}\|_1 + o(1)\\
&\le 3\int_{[f\ne g]}(1+\phi\circ T_{-n})d\mu+ o(1)\\ &\le
3(\d+\e)+o(1).\endalign$$ Thus $$\|g\circ T_n\circ\pi\circ
T_{-n}-g\|_1=\|g\circ T_n\circ\pi-g\circ T_{n}\|_1\to 0$$
establishing ($\star$).\f a) $\Leftarrow$) Suppose that
$\nu\in\Cal P_e(\Cal P_e(X,\goth T(X),\Psi_G),S)$. Evidently
$$\mu_\nu\circ T_k=\int_{\Cal P_e(X,\goth T(X),\Psi_G)}\om\circ T_k
d\nu(\om)=\int_{\Cal P_e(X,\goth T(X),\Psi_G)}S_k\om
d\nu(\om)=\mu_\nu.$$ To check $T$-ergodicity let $A\in\B(X),\
T_kA=A\ \forall\ k\in\Bbb Z^d$. Let $\pi\in [\goth T(X)]$. It
follows from ($\star$) that $$0\ \leftarrow\ \mu_\nu(T_n\pi
T_{-n}A\D A)=\mu_\nu(\pi A\D A)$$ and  $1_A\circ \pi=1_A\
\mu_\nu$-a.e., whence $\om$-a.e. for $\nu$-a.e. $\om\in\Cal
P_e(X,\goth T(X),\Psi_G)$. Thus $\om(A)=0,1$ for $\nu$-a.e.
$\om\in\Cal P_e(X,\goth T(X),\Psi_G)$. The set $\tilde
A:=\{\om\in\Cal P_e(X,\goth T(X),\Psi_G):\ \om(A)=1\}$ is
$S$-invariant, whence $\nu(\tilde A)=0,1$ and $\mu_\nu(A)=0,1$.
 \smallskip\f a) $\Rightarrow$)

 Let $\mu=\mu_\nu\in\Cal P(X)$ be a $G$-equilibrium measure
where $\nu\in\Cal P(\Cal P_e(X,\goth T(X),\Psi_G))$. As above,
$\nu$ is $S$-invariant. The $S$-ergodicity on $\nu$ follows from
the identification of Borel measurable $S$-invariant subsets of
$\Cal P_e(X,\goth T(X),\Psi_G)$ with $\goth I(\goth
T(X))$-measurable $T$-invariant subsets of $X$.
\par b) By ($\star$), $$\mu(T_{-n}\pi T_nA\D
A)\to 0\ \text{ as }\ n\to\infty\ \forall\ A\in\B(X),\ \pi\in
[\goth T(X)].$$

Now let $A\in\B(X)$ be a $T$-rigid set. We show that $A\in\goth
I(\goth T(X))\ \mod\ \mu$. To this end  suppose that
$n_k\to\infty,\ \mu(T_{n_k}A\D A)\to 0$, and let  $\pi\in [\goth
T(X)]$, then
$$0\ \leftarrow\ \mu(T_{-n_k}\pi T_{n_k}A\D A)=\mu(\pi T_{n_k}A\D
T_{n_k}A)\to \mu(\pi A\D A),$$ $\pi A=A\ \mod\ \mu$ and $A\in\goth
I(\goth T(X))\ \mod\ \mu$ which is identified with the factor
$\s$-algebra.
\par c) It follows from b) that any $T$-rigid set is the pull-back
of an $S$-rigid set, and thus from the remarks preceding this
proposition that each $T$-eigenfunction is the pull-back of a
$S$-eigenfunction with the same eigenvalue. Statement c) follows
from this.
 \hfill\qed\enddemo
 \subheading{The Burton Steif construction}

  Let $X\subset S^{\Bbb Z^d}$
be a subshift. As in \cite{BS2}, define the {\it free $\Bbb Z$-product
of $X$} by
$$Z:=\{x\in S^{\Bbb Z^{d+1}}:\ x^{(n)}\in X\ \forall\ n\in\Bbb Z\}$$
where for $x\in S^{\Bbb Z^{d+1}},\ n\in\Bbb Z$, $x^{(n)}\in
S^{\Bbb Z^{d}}$ is defined by $x^{(n)}_k:=x_{(k,n)}\ \ (k\in\Bbb
Z^d)$. Evidently $Z$ is a $\Bbb Z^{d+1}$-subshift, and if $X$ is
strongly irreducible and/or strongly aperiodic TMS, then so is $Z$.

 It was shown in
\cite{BS2} that if $X$ has more than one equilibrium measure, then
$Z$ has uncountably many.

Here, we study the collection of equilibrium measures for $Z$
using proposition 4.1. The tail relation $\goth T(Z)$ has a
product structure. \subheading{Product relations}
 Suppose that $\r$ is a countable
equivalence relation on the Polish space $Y$ and that $\G$ is an
at most countable set. The {\it $\G$-product of \ \ $\r$} is the
equivalence relation $\r^{(\G)}\in\B(Y^\G\x Y^\G)$ defined by
$$\align &\r^{(\G)}:=\\ &
\{(x,x')\in Y^\G\x Y^\G:\ \exists\ F\subset\G\ \text{finite},\
(x_\g,x'_\g)\in\r\ \forall\ \g\in F,\
x_{F^c}=x'_{F^c}\}.\endalign$$ \proclaim{Lemma 4.2} If  $X$ is a
TMS and $Z$ is its free $\Bbb Z$-product, then $Z=X^\Bbb Z$ and
$\goth T(Z)=\goth T(X)^{(\Bbb Z)}.$ \endproclaim  Now suppose that
$\Psi:\r\to\Bbb R_+$ is a multiplicative $\r$-cocycle. Define the
{\it product cocycle} $\Psi^{(\G)}:\r^{(\G)}\to\Bbb R_+$ by
$$\Psi^{(\G)}(x,y):=\prod_{\g\in\G}\Psi(x_\g,y_\g).$$
This formula defines a multiplicative $\r^{(\G)}$-cocycle as each
product only has finitely many non-unit terms.

Let $\mu_\g\in\Cal P(Y)\ \ (\g\in\G)$ and set
$\mu:=\prod_{\g\in\G}\mu_\g\ \in\Cal P(Y^\G)$.
 \proclaim{Lemma 4.3}$$\mu\in\Cal P_e(Y^\G,\r^{(\G)},\Psi^{(\G)})\
 \iff\ \mu=\prod_{\g\in\G}\mu_\g\
\text{\rm where}\ \mu_\g\in\Cal P_e(Y,\r,\Psi)\ \forall\
\g\in\G.$$
\endproclaim\demo{Proof}

Let $\mu_\g\in\Cal P(Y)\ \ (\g\in\G)$ and set
$\mu:=\prod_{\g\in\G}\mu_\g\ \in\Cal P(Y^\G)$. It is routine to
show that
$$\mu\in\Cal P(Y^\G,\r^{(\G)},\Psi^{(\G)})\
 \iff\ \mu_\g\in\Cal P(Y,\r,\Psi)\ \forall\
\g\in\G.$$ We turn to the ergodicity assertions. \smallskip\f
$\Leftarrow$)

Suppose that $\mu=\prod_{\g\in\G}\mu_\g$ where $\mu_\g\in\Cal
P_e(Y,\r,\Psi)\  (\g\in\G)$ and let $A\in\B(Y^\G)$ be
$\r^{(\G)}$-invariant. \par For $\g_0\in\G$ and $x\in
Y^{\G\setminus\{\g_0\}}$ let $A_{\g_0,x}:=\{y\in Y:\ (x,y)\in
A\}$, then  $A_{\g_0,x}$ is $\r$-invariant $\forall\ x\in
Y^{\G\setminus\{\g_0\}}$.

Since $\mu_{\g_0}$ is $\r$-ergodic, $\mu_{\g_0}(A_{\g_0,x})=0,1$
for $\prod_{\g\in\G\setminus\{\g_0\}}\mu_\g$-a.e. $x\in
Y^{\G\setminus\{\g_0\}}$, and $A\overset{\mu}\to{=}\{x\in
Y^{\G\setminus\{\g_0\}}:\
\mu_{\g_0}(A_{\g_0,x})=1\}\in\B(Y^\G\setminus\{\g_0\})$.
\par Continuing analogously shows that for any finite set
$F\subset\G,\ A\in\B(Y^{\G\setminus F})$. But then $A$ is
$\mu$-independent of every set $B\in\s\(\bigcup_{F\subset\G\
\text{\rm finite}}\B(Y^F)\)=\B(Y^\G)$ and $\mu(A)=0,1$.
\smallskip\f $\Rightarrow$)

 Suppose that $\mu\in\Cal P_e(Y^\G,\r^{(\G)},\Psi^{(\G)})$. Note that for any $E\subset\G$,
 $$\mu_E:=\mu\circ\pi_E^{-1}\in\Cal P_e(Y^E,\r^{(E)},\Psi^{(E)})$$
 where $\pi_E(x):=x_E$.
 We'll show that $\mu=\mu_E\x\mu_{\G\setminus E}\ \forall\ E\subset\G$.
 \par Denoting $Y^\G=Y^E\x Y^{\G\setminus E}$, we have by the disintegration
 theorem that
  $$\mu(A\x B)=\int_A\nu_x(B)d\mu_E(x)\ \ \ (A\in\B(Y^E),\ B\in\B(Y^{\G\setminus
  E}))$$ where $x\mapsto\nu_x$ is a measurable mapping ($Y^E\to\Cal
  P(Y^{\G\setminus
  E})$). Let $A_0\in \B(Y^E)$ and let $V:A_0\to VA_0$ be a
  $\r^{(E)}$-holonomy. Let $\widetilde V:A_0\x Y^{\G\setminus
  E}\to V(A_0)\x Y^{\G\setminus
  E}$ be the corresponding $\r^{(\G)}$-holonomy defined by
  $\widetilde V(x_E,x_{G\setminus
  E}):=(V(x_E),x_{G\setminus
  E})$, then
  $$\tfrac{d\mu\circ\widetilde V}{d\mu}(x)=\Psi^{(\G)}(x,\widetilde Vx)=
  \Psi^{(E)}(x,Vx_E)=\tfrac{d\mu_E\circ V}{d\mu_E}(x_E).$$

  Thus,   $\forall\ A\in\B(Y^E),\ A\subset
  A_0$ and $B\in\B(Y^{\G\setminus
  E})$,
  $$\align \int_A\nu_{V^{-1}x}(B)d\mu_E(x) &=
  \int_{X_E}\tfrac{d\mu_E\circ V}{d\mu_E}(x)1_A\circ V(x)\nu_{x}(B)d\mu_E(x)\\ &=
  \int_X(\tfrac{d\mu_E\circ V}{d\mu_E}1_A\circ V)\otimes 1_Bd\mu\\ &=
\int_X\tfrac{d\mu\circ\widetilde V}{d\mu}1_{A\x B}\circ\widetilde
V d\mu\\ &=
   \mu(A\x
  B)\\ &=\int_{A}\nu_x(B)d\mu_E(x)\endalign$$
  and $x\mapsto\nu_x(B)$ is $\r^{(E)}$-invariant, whence
  $\mu_E$-a.e. constant $\forall\ B\in\B(Y^{\G\setminus
  E})$. It follows that for $\mu_E$-a.e. $x\in X_E,$
  $$\nu_x=\int_{X_E}\nu_yd\mu_E(y)=\mu_{\G\setminus E}.$$\hfill\qed\enddemo

\demo{ Proof of theorem 3}
 \par For 1), fix  a strongly irreducible, strongly aperiodic $\Bbb Z^2$-TMS $X$,
 for which  there are
 two equilibrium measures $P_+,P_-$ which are $T$-weakly
Bernoulli (e.g. suitable iceberg, or beach models as in \S3, see \cite{BS2}). It follows from \cite{Ho-St1} that $P_+,P_-$ are both
$\goth T(X)$-ergodic.

For 2), let $G:X:=\{-1,1\}^{\Bbb Z^2}\to\Bbb R$ be the Markov
potential defined by $G(x):=\b x_{n}x_{n+e_1}x_{n+e_2}$ where
$\b>0$. As shown in \cite{P} for $\b$ large enough, there are
 two $G$-equilibrium measures $P_+,P_-$ which are shown to be $T$-weakly
Bernoulli in \cite{LGR}, whence $\goth T(X)$-ergodic by
\cite{Ho-St2}. Fix such $\b>0$.

In both cases, let $Z$ be the free $\Bbb Z$-product of $X$, a
$\Bbb Z^3$-TMS, define $\widetilde G:Z\to\Bbb R$ by $$\widetilde
G(x)=\cases & 0\ \ \ \ \text{\rm in case 1)},\\ & \b
x_{n}x_{n+e_1}x_{n+e_2}\ \ \ \ \text{\rm in case 2)}.\endcases$$

\par By lemma 4.3,
 for each $\eta\in\{-,+\}^{\Bbb Z}$,
 $$P_{\eta}:=\prod_{\ell\in\Bbb Z}P_{\eta_\ell}\in\Cal P_e(Z,\goth T(Z),e^{\Psi_{\widetilde G}}).$$
 Also, $P_{\eta}\circ T_k=P_{\eta}\ \forall\
k\in\Bbb Z^2\x\{0\}$. Thus, $\{-,+\}^{\Bbb Z}$ is invariant under each $S_{(k_1,k_2,k_3)}$ and
$S_{(k_1,k_2,k_3)}|_{\{-,+\}^{\Bbb
Z}}=\s^{k_3}$ where $\s:\{-,+\}^{\Bbb Z}\to\{-,+\}^{\Bbb Z}$ is
the shift.\par   If $\nu\in\Cal P(\{-,+\}^{\Bbb Z})$ then
$P_\nu:=\int_{\{-,+\}^{\Bbb Z}}P_{\eta}d\nu(\eta)\in\Cal P(Z,\goth
T(Z),e^{\Psi_{\widetilde G}})$. By proposition 4.1, $P_\nu$ is a
$\widetilde G$-equilibrium measure iff $\nu$ is $\s$-invariant,
ergodic.

Now let $\rho\in\Cal P(\Bbb T)$ be the spectral type of some
ergodic, probability preserving $\Bbb Z$-action. As is well known,
$\exists\ \nu\in\Cal P_e(\{-,+\}^{\Bbb Z},\s)$ so that
$(\{-,+\}^{\Bbb Z},\nu,\s)$ has spectral type $\rho$ (i.e. any
spectral type can be achieved by an ergodic, probability
preserving $\Bbb Z$-action with entropy less than $\log 2$).

 We complete the proof of spectral  abundance by showing how the spectral properties of the
dynamical system $(Z,T,P_\nu)$ reflect those of $(\{-,+\}^{\Bbb
Z},\s,\nu)$.
\par For any $\eta\in\{-,+\}^{\Bbb Z}$, the $\Bbb Z^2$-random
field $(Z,T|_{\Bbb Z^2\x\{0\}},P_\eta)$ is weakly Bernoulli,
whence strongly mixing. Thus for $\nu\in\Cal P_e(\{-,+\}^{\Bbb
Z},\s)$,
$$P_\nu(A\cap
T_{(n_1,n_2,0)}B)\underset{(n_1,n_2)\to\infty}\to\longrightarrow
\int_{\{-,+\}^{\Bbb Z}}P_{\eta}(A)P_{\eta}(B)d\nu(\eta)\ \forall\
A,B\in\B(Z).\tag1$$ \par The measures $\prod P_+,\ \prod
P_-\in\Cal P_e(Z,\goth T(Z),e^{\Psi_{\widetilde G}})$ are weak
Bernoulli equilibrium measures on $Z$. \par If $\nu\in\Cal
P_e(\{-,+\}^{\Bbb Z},\s),\ \nu\ne\d_{(+,+,\dots)},\
\d_{(-,-,\dots)}$, then $(Z,T,P_\nu)$ is not strongly mixing
because: $\nu$ is not a point mass whence $\int_{\{-,+\}^{\Bbb
Z}}P_{\eta}(A)^2d\nu(\eta)> P_\nu(A)^2$ whenever $0<P_\nu(A)<1$
and strong mixing of $(Z,P_\nu,T)$ is eliminated by (1).
\par By part b) of proposition 4.1, the weak mixing, mild mixing, total ergodicity of $T$ is equivalent to that of
$S$, which in turn is equivalent to that of $\s$ (respectively).
\hfill\qed\enddemo

  \heading References\endheading

\enddocument